\PassOptionsToPackage{unicode=true}{hyperref} % options for packages loaded elsewhere
\PassOptionsToPackage{hyphens}{url}
\documentclass[]{article}
\usepackage{lmodern}
\usepackage{amssymb,amsmath}
\usepackage{ifxetex,ifluatex}
\usepackage{fixltx2e} % provides \textsubscript
\ifnum 0\ifxetex 1\fi\ifluatex 1\fi=0 % if pdftex
  \usepackage[T1]{fontenc}
  \usepackage[utf8]{inputenc}
  \usepackage{textcomp} % provides euro and other symbols
\else % if luatex or xelatex
  \usepackage{unicode-math}
  \defaultfontfeatures{Ligatures=TeX,Scale=MatchLowercase}
\fi
\IfFileExists{upquote.sty}{\usepackage{upquote}}{}
\IfFileExists{microtype.sty}{%
\usepackage[]{microtype}
\UseMicrotypeSet[protrusion]{basicmath} % disable protrusion for tt fonts
}{}
\IfFileExists{parskip.sty}{%
\usepackage{parskip}
}{% else
\setlength{\parindent}{0pt}
\setlength{\parskip}{6pt plus 2pt minus 1pt}
}
\usepackage{hyperref}
\hypersetup{
            pdfborder={0 0 0},
            breaklinks=true}
\usepackage{longtable,booktabs}
\IfFileExists{footnote.sty}{\usepackage{footnote}\makesavenoteenv{longtable}}{}
\ifx\paragraph\undefined\else
\let\oldparagraph\paragraph
\renewcommand{\paragraph}[1]{\oldparagraph{#1}\mbox{}}
\fi
\ifx\subparagraph\undefined\else
\let\oldsubparagraph\subparagraph
\renewcommand{\subparagraph}[1]{\oldsubparagraph{#1}\mbox{}}
\fi

\makeatletter
\def\fps@figure{htbp}
\makeatother

\date{}

\begin{document}

\textbf{A Unified Exact Factorial-Moment Theory for Multi-set Allocation
Occupancy (MAO) in Finite Populations}

Xing-gang Mao\textsuperscript{a*} and Xiao-yan Xue\textsuperscript{b}

\textsuperscript{a}Department of Neurosurgery, Xijing Hospital, the
Fourth Military Medical University, Xi'an, No. 17 Changle West Road,
Xi'an, Shaanxi Province, China;

\textsuperscript{b}Department of Pharmacology, School of Pharmacy, the
Fourth Military Medical University, No. 17 Changle West Road, Xi'an,
Shaanxi Province, China

*Corresponding authors. E-mail addresses:
{\href{mailto:xgmao@fmmu.edu.cn}{\nolinkurl{xgmao@fmmu.edu.cn}} (X.G.
Mao) }

Author ORCID: X.G. Mao: 0000-0001-8505-1904.

\textbf{\\
}

\textbf{Abstract}

Let~\(A_{1},\ldots,A_{T}\)~be independent uniformly selected subsets of
a finite population of size~\(n\), with prescribed
cardinalities~\(m_{1},\ldots,m_{T}\). For each population element,
define its occupancy level as the number of selected subsets containing
it. Let~\(x_{t}\)~and~\(x_{\geq t}\)~denote the numbers of elements with
occupancy exactly~\(t\)~and at least~\(t\), respectively.

The 2025 work introduced a general multi-set allocation occupancy (MAO)
representation for higher-order occupancy moments. The present paper
establishes its joint-probabilistic interpretation and provides a
rigorous unified derivation for arbitrary joint occupancy categories and
moment orders. Specifically, for
arbitrary~\(B_{1},\ldots,B_{\mathcal{l}} \subseteq \left\{ 0,1,\ldots,T \right\}\),
we prove the exact representation

\[F_{\mathcal{l}}\left( B_{1},\ldots,B_{\mathcal{l}} \right) = \frac{G_{T}\left( B_{1},\ldots,B_{\mathcal{l}} \right)}{\left( n \right)_{\mathcal{l}}^{T - 1}},\]

where~\(G_{T}\left( B_{1},\ldots,B_{\mathcal{l}} \right)\)~is the
corresponding generalized MAO transversal sum. This identity unifies the
joint factorial moments of arbitrary occupancy categories within a
single exact finite-population framework.

In particular,
writing~\(B_{\geq t} = \left\{ t,t, + ,1\ldots T \right\}\), the
factorial moments of exact- and threshold-occupancy counts are obtained
as the specializations

\[\mathbb{E}\left\lbrack \left( x_{t} \right)_{\mathcal{l}} \right\rbrack = F_{\mathcal{l}}\left( \left\{ t \right\},\ldots,\left\{ t \right\} \right),\]

and

\[\mathbb{E}\left\lbrack \left( x_{\geq t} \right)_{\mathcal{l}} \right\rbrack = F_{\mathcal{l}}\left( B_{\geq t},\ldots,B_{\geq t} \right).\]

Mixed factorial moments, raw moments, variances, and covariances follow
from the same representation through standard transformations. The
formulas are verified by exhaustive enumeration over feasible parameter
ranges and by Monte Carlo simulation. The resulting theory provides a
rigorous and unified finite-population foundation for exact and
threshold multi-set occupancy statistics.


\textbf{1. Introduction}

Repeated membership across multiple subsets of a common finite
population is a basic structure in combinatorial probability, occupancy
theory, and finite-population statistics. Let several subsets of
prescribed sizes be selected from a common universe. The resulting basic
question is: how many elements occur in exactly~\(t\), or in at
least~\(t\), of the selected subsets? This structure occurs in repeated
feature selection, replicated experiments, screening procedures,
independent cohort studies, network analysis, ensemble methods, and
multi-set overlap analysis.

The problem is nevertheless mathematically difficult. Fixed-size subset
selection is sampling without replacement, and membership indicators of
distinct population elements are therefore dependent.
With~\(T\)~subsets, each element may have any of~\(2^{T}\)~membership
patterns. Exact analysis must simultaneously respect the
finite-population constraint, prescribed subset sizes, dependence across
elements, and aggregation over many membership patterns. These are
characteristic difficulties of finite-population and combinatorial
probability (Cochran, 1991; Hájek, 1960; Särndal et al., 1992).

Let

\[\Omega = \{ 1,\ldots,n\}\]

be a finite population, and let

\[A_{1},\ldots,A_{T} \subseteq \Omega\]

be mutually independent random subsets, where~\(A_{i}\)~is uniformly
distributed over all~\(m_{i}\)-element subsets of~\(\Omega\). Thus,

\[\mid A_{i} \mid = m_{i},\ 1 \leq i \leq T.\]

For~\(a \in \Omega\), define the occupancy level

\[L\left( a \right) = \sum_{i = 1}^{T}1_{\left\{ a \in A_{i} \right\}}.\]

The exact-occupancy and threshold-occupancy counts are

\[x_{t} = \#\left\{ a \in \Omega:L\left( a \right) = t \right\},\ x_{\geq t} = \#\left\{ a \in \Omega:L\left( a \right) \geq t \right\}.\]

For~\(T = 2\), the fully occupied
count~\(x_{2} = \mid A_{1} \cap A_{2} \mid\)~has the classical
hypergeometric distribution. For~\(T \geq 3\), pairwise intersections no
longer determine the allocation of elements among all subsets, and the
distributions of~\(x_{t}\)~and~\(x_{\geq t}\)~become genuinely multi-set
finite-population problems.

This model was originally introduced as the~\textbf{general
hypergeometric distribution}~(GHGD) (Mao and Xue, 2018; Mao and Xue,
2022), and named as \textbf{Multi-set Allocation Occupancy (MAO)
distribution}~in the following works (Mao, 2026; Mao and Xue, 2025).
GHGD emphasizes its extension of the classical hypergeometric law, while
MAO emphasizes the underlying allocation of each population element to a
membership pattern and its resulting occupancy level. In addition, there
are works studied the \(x_{0}\) problems (Vatutin and Mikhajlov, 1982),
which is a special situation of the MAO model (distribution of
\(x_{0}\)).

The central result of this paper is giving unified exact closed-form
formulas for any order moments of the distribution of
\(x_{t}\)~and~\(x_{\geq t}\), and providing a strict mathematical proof
based on statistical induction. For arbitrary value of subset count
\(T\), occupancy level~\(t\), prescribed size
vector~\(\left( m_{1},\ldots,m_{T} \right)\), and moment order~\(r\),
the factorial moments of both

\[x_{t}\ \mathrm{\text{and}}\ x_{\geq t}\]

are given by one common MAO factorial moments.

More precisely, we prove that MAO allocation weights are unnormalized
joint probabilities of prescribed membership patterns for distinct
population elements. Their normalized sums over occupancy-compatible
patterns are factorial moments. The MAO norms introduced previously
(Mao, 2026; Mao and Xue, 2025) are therefore precisely MAO factorial
moments. Ordinary moments follow through the Stirling transformation

\[y^{v} = \sum_{r = 0}^{v}\begin{Bmatrix}
v \\
r \\
\end{Bmatrix}\left( y \right)_{r}.\]

Hence the same formula yields exact means, variances, covariances, any
order of raw moments, and central moments.
Since~\(x_{t}\)~and~\(x_{\geq t}\)~have finite support, their complete
moment sequences determine their distributions.

The present paper consolidates and rigorously completes a sequence of
developments in the GHGD/MAO framework. The 2018 work introduced the
generalized hypergeometric distribution and recursive enumeration
procedures. The 2022 work established exact formulas for the means and
variances of the occupancy counts for arbitrary~\(T\)~and~\(t\). The
2025 work introduced the MAO allocation weights and MAO norms, thereby
exposing a compact higher-order algebraic structure. The present paper
proves the general joint-probabilistic interpretation of this structure
and establishes the corresponding arbitrary-order factorial-moment
representation. In this way, the earlier distributional, moment, and MAO
constructions are unified into a single rigorous finite-population
theory.

The present paper consolidates and rigorously completes a sequence of
developments in the GHGD/MAO framework (Table 1). The 2018 work
introduced the generalized hypergeometric distribution and recursive
enumeration procedures (Mao and Xue, 2018) . The 2022 work established
exact formulas for the means and variances of the occupancy counts for
arbitrary~\(T\)~and~\(t\) (Mao and Xue, 2022). The 2025 work introduced
the MAO allocation weights and MAO norms, thereby exposing a compact
higher-order algebraic structure (Mao and Xue, 2025). The present paper
proves the general joint-probabilistic interpretation of this structure
and establishes the corresponding arbitrary-order factorial-moment
representation. In this way, the earlier distributional, moment, and MAO
constructions are unified into a single rigorous finite-population
theory.

\textbf{Table 1. Development of the GHGD/MAO theory.}

\begin{longtable}{@{}p{0.22\textwidth}p{0.36\textwidth}p{0.42\textwidth}@{}}
\caption{Development of the GHGD/MAO theory.}
\label{tab:development}\\
\toprule
\textbf{Work} &
\textbf{Main result} &
\textbf{Limitation resolved in the present paper} \\
\midrule
\endfirsthead

\toprule
\textbf{Work} &
\textbf{Main result} &
\textbf{Limitation resolved in the present paper} \\
\midrule
\endhead

\midrule
\multicolumn{3}{r@{}}{Continued on the next page} \\
\endfoot

\bottomrule
\endlastfoot

2018 (Mao and Xue, 2018) &
Introduced GHGD and recursive counting; derived initial moment results. &
No uniform all-order formula. \\

2022 (Mao and Xue, 2022) &
Exact means and variances for arbitrary $T$ and $t$. &
No common arbitrary-order moment mechanism. \\

2025 (Mao and Xue, 2025) &
Introduced MAO allocation weights and MAO norms. &
Probability interpretation and rigorous general derivation unresolved. \\

Present paper &
Exact MAO factorial moments of arbitrary order. &
Consolidates the GHGD/MAO framework into a unified rigorous moment theory. \\

\end{longtable}

All closed-form formulas derived here were checked against exact
exhaustive enumeration whenever enumeration is computationally feasible,
with exact agreement. In larger parameter regimes, the theoretical
results were additionally verified by Monte Carlo simulation, showing
close agreement within simulation error.

One important application of the MAO distribution is the statistical
interpretation of multi-set overlaps commonly represented by Venn
diagrams and related visualizations (Amand et al., 2019; Hur et al.,
2019). Such diagrams are widely used, especially when repeated
experiments seek elements reproduced in exactly~\(t\), or at
least~\(t\), among~\(T\)~studies. They display observed overlap patterns
(Lex et al., 2014), but do not provide null distributions for these
repeated-occurrence counts. The MAO distribution supplies this missing
inferential foundation.

Accordingly, the MAO distribution has potential applications to a broad
range of finite-population, multi-set overlap problems requiring
statistical inference. Relevant areas include combinatorics and extremal
set systems (Alon and Spencer, 2016), network science and
overlapping-community analysis (Fortunato, 2010), genomics,
transcriptomics, and multi-study biomarker replication (Pavlopoulos et
al., 2018), as well as ecology and environmental science, especially
species co-occurrence and niche-overlap analysis (Wilson and Anderson,
1985). Further applications may arise in neuroscience and psychology
(Borsboom et al., 2021), where reproducibility is assessed across
studies or analytic pipelines; in social, economic, and bibliometric
research involving repeated participation, affiliation, or
co-occurrence; and in physics and materials science, where the recurrent
spatial occupancy of particles, defects, or other features is examined
across repeated observations (Dokuchaev, 2000).

\textbf{2. The MAO Sampling Model}

Let~\(\Omega\)~be a finite population with

\[\mid \Omega \mid = n.\]

Let~\(T \geq 2\), and for
each~\(i \in \left\lbrack T \right\rbrack: = \{ 1,\ldots,T\}\), let

\[A_{i} \subseteq \Omega,\  \mid A_{i} \mid = m_{i},\ 0 \leq m_{i} \leq n.\]

\[\mathbf{m}: = \left( m_{1},\ldots,m_{T} \right)\]

The random sets~\(A_{1},\ldots,A_{T}\)~are independently sampled, with
each~\(A_{i}\)~uniformly distributed over
all~\(\left( \frac{n}{m_{i}} \right)\)~subsets of~\(\Omega\)~of
size~\(m_{i}\).

Thus, the sample space contains

\[\prod_{i = 1}^{T}\left( \frac{n}{m_{i}} \right)\]

ordered configurations, each having equal probability.

For each~\(a \in \Omega\), define its occupancy level by

\[L\left( a \right): = \sum_{i = 1}^{T}1_{\left\{ a \in A_{i} \right\}}.\]

Hence,

\[L\left( a \right) \in \{ 0,1,\ldots,T\}.\]

\textbf{Definition 1. Exact and thresholded MAO counts}

For~\(t \in \{ 0,\ldots,T\}\), define

\[x_{t}: = \sum_{a \in \Omega}^{}1_{\left\{ L\left( a \right) = t \right\}}.\]

Thus,~\(x_{t}\)~is the number of elements occurring in
exactly~\(t\)~selected subsets.

Define

\[x_{\geq t}: = \sum_{a \in \Omega}^{}1_{\left\{ L\left( a \right) \geq t \right\}} = \sum_{q = t}^{T}x_{q}.\]

The MAO distribution refers to the joint distribution of the occupancy
vector

\[\left( x_{0},x_{1},\ldots,x_{T} \right),\]

as well as to the marginal distributions
of~\(x_{t}\)~and~\(x_{\geq t}\).

For every realization,

\[\sum_{t = 0}^{T}x_{t} = n,\]

and

\[\sum_{t = 0}^{T}tx_{t} = \sum_{i = 1}^{T}m_{i}.\]

When~\(T = 2\),~\(x_{2} = \mid A_{1} \cap A_{2} \mid\)~follows the
classical hypergeometric distribution:

\[x_{2} \sim \text{Hypergeometric}\left( n,m_{1},m_{2} \right).\]

\textbf{3. One-Element Occupancy Distribution}

Fix~\(a \in \Omega\), and define the membership indicators

\[I_{i}\left( a \right): = 1_{\left\{ a \in A_{i} \right\}},\ i \in \left\lbrack T \right\rbrack.\]

Since~\(A_{i}\)~is uniformly sampled from the subsets of~\(\Omega\)~of
size~\(m_{i}\),

\[\mathbb{P}\left\{ I_{i}\left( a \right) = 1 \right\}\mathbb{= P}\left\{ a \in A_{i} \right\} = \frac{m_{i}}{n}\mathbb{,\ P}\left\{ I_{i}\left( a \right) = 0 \right\} = \frac{n - m_{i}}{n}.\]

Hence,

\[I_{i}\left( a \right) \sim \text{Bernoulli}\mathrm{\,}\left( \frac{m_{i}}{n} \right).\]

Its probability generating function (PGF) is

\[Q_{i}\left( z \right)\mathbb{: = E}\mathrm{\,}\left\lbrack z^{I_{i}\left( a \right)} \right\rbrack = \frac{n - m_{i}}{n} + \frac{m_{i}}{n}\text{z.}\]

Because~\(A_{1},\ldots,A_{T}\)~are independent, the random variables
\(I_{1}\left( a \right),\ldots,I_{T}\left( a \right)\)~are independent.
Therefore, for

\[L\left( a \right): = \sum_{i = 1}^{T}I_{i}\left( a \right),\]

the PGF of~\(L\left( a \right)\)~is

\[\begin{matrix}
 & \begin{matrix}
Q\left( z \right)\mathbb{\&: = E}\mathrm{\,}\left\lbrack z^{L\left( a \right)} \right\rbrack \\
\& = \prod_{i = 1}^{T}\mathbb{E}\mathrm{\,}\left\lbrack z^{I_{i}\left( a \right)} \right\rbrack \\
\& = \prod_{i = 1}^{T}\left( \frac{n - m_{i}}{n} + \frac{m_{i}}{n}z \right). \\
\end{matrix} & & \mathrm{(1)} \\
\end{matrix}\]

For a polynomial~\(f\left( z \right)\),
let~\(\left\lbrack z^{t} \right\rbrack f\left( z \right)\)~denote the
coefficient of \(z^{t}\). Then

\[\mathbb{P}\left\{ L\left( a \right) = t \right\} = \left\lbrack z^{t} \right\rbrack Q\left( z \right),\ t = 0,\ldots,T.\]

\textbf{Theorem 1. Exact MAO means}

For every~\(t \in \{ 0,\ldots,T\}\),

\[\begin{matrix}
 & \mathbb{E}\left\lbrack x_{t} \right\rbrack = n\left\lbrack z^{t} \right\rbrack Q\left( z \right). & & \mathrm{(2)} \\
\end{matrix}\]

Equivalently,

\[\begin{matrix}
 & \mathbb{E}\left\lbrack x_{t} \right\rbrack = n\sum_{\substack{S \subseteq \left\lbrack T \right\rbrack \\ \mid S \mid = t}}^{}{\left( \prod_{i \in S}^{}\frac{m_{i}}{n} \right)\left( \prod_{i \notin S}^{}\frac{n - m_{i}}{n} \right)}. & & \mathrm{(3)} \\
\end{matrix}\]

\textbf{Proof.}~Define

\[I_{a}^{\left( t \right)}: = 1_{\left\{ L\left( a \right) = t \right\}}.\]

Then

\[x_{t} = \sum_{a \in \Omega}^{}I_{a}^{\left( t \right)}.\]

Therefore, by linearity of expectation,

\[\mathbb{E}\left\lbrack x_{t} \right\rbrack = \sum_{a \in \Omega}^{}\mathbb{E}\left\lbrack I_{a}^{\left( t \right)} \right\rbrack = \sum_{a \in \Omega}^{}\mathbb{P}\left( L\left( a \right) = t \right).\]

The MAO model is invariant under permutations of the~\(n\)~population
elements, so the distribution of~\(L\left( a \right)\)~is the same for
every~\(a \in \Omega\). Hence,

\[\mathbb{E}\left\lbrack x_{t} \right\rbrack = n\mathrm{\,}\mathbb{P}\left( L\left( a \right) = t \right).\]

Because~\(Q\left( z \right)\)~is the PGF of~\(L\left( a \right)\),

\[\mathbb{P}\left( L\left( a \right) = t \right) = \left\lbrack z^{t} \right\rbrack Q\left( z \right).\]

This proves (2). Expanding the product in (1) and collecting the
coefficient of~\(z^{t}\)~yields (3).~~

\textbf{Corollary 1. Mean of the fully occupied count}

\[\begin{matrix}
 & \mathbb{E}\left\lbrack x_{T} \right\rbrack = \frac{\prod_{i = 1}^{T}m_{i}}{n^{T - 1}}. & & \mathrm{(4)} \\
\end{matrix}\]

\textbf{Corollary 2. Mean of the thresholded count}

For every~\(t \in \{ 0,\ldots,T\}\),

\[\begin{matrix}
 & \mathbb{E}\left\lbrack x_{\geq t} \right\rbrack = n\sum_{q = t}^{T}\lbrack\left. \ z^{q} \right\rbrack Q\left( z \right). & & \mathrm{(5)} \\
\end{matrix}\]

\textbf{Corollary 3. Equal-size subsets}

If

\[m_{1} = \cdots = m_{T} = m,\]

then

\[Q\left( z \right) = \left( 1 - \frac{m}{n} + \frac{m}{n}z \right)^{T},\]

and

\[\begin{matrix}
 & \mathbb{E}\left\lbrack x_{t} \right\rbrack = n\left( \frac{T}{t} \right)\left( \frac{m}{n} \right)^{t}\left( 1 - \frac{m}{n} \right)^{T - t}. & & \mathrm{(6)} \\
\end{matrix}\]

\textbf{4. Two-Element Occupancy Distribution and Exact Second Moments}

For~\(t \in \{ 0,\ldots,T\}\), recall that

\[I_{a}^{\left( t \right)}: = 1_{\left\{ L\left( a \right) = t \right\}},\ x_{t} = \sum_{a \in \Omega}^{}I_{a}^{\left( t \right)}.\]

Although the occupancy indicators associated with a fixed element are
determined by its membership across the independently sampled sets, the
indicators associated with two distinct population elements are
generally dependent. More precisely, for~\(a \neq b\), the
variables~\(I_{a}^{\left( r \right)}\)~and~\(I_{b}^{\left( s \right)}\)~are
generally not independent. The dependence arises from the fixed capacity
of each sampled set~\(A_{i}\): conditional on one specified element
being included in~\(A_{i}\), fewer of the~\(m_{i}\)~available positions
remain for the other specified element.

Fix distinct elements

\[a,b \in \Omega.\]

For a single uniformly sampled subset~\(A_{i} \subseteq \Omega\)~of
cardinality~\(m_{i}\), the four possible joint membership states
of~\(a\)~and~\(b\)~have probabilities

\[\begin{matrix}
 & P(a \notin A_{i},b \notin A_{i}) = \frac{\left( n - m_{i} \right)\left( n - m_{i} - 1 \right)}{n\left( n - 1 \right)}, & & \mathrm{(}\mathrm{7}\mathrm{)} \\
\end{matrix}\]

\[\begin{matrix}
 & P(a \in A_{i},b \notin A_{i}) = P(a \notin A_{i},b \in A_{i}) = \frac{m_{i}\left( n - m_{i} \right)}{n\left( n - 1 \right)}, & & \mathrm{(}\mathrm{8}\mathrm{)} \\
\end{matrix}\]

and

\[\begin{matrix}
 & P(a \in A_{i},b \in A_{i}) = \frac{m_{i}\left( m_{i} - 1 \right)}{n\left( n - 1 \right)}. & & \mathrm{(}\mathrm{9}\mathrm{)} \\
\end{matrix}\]

Define the per-set joint PGF

\[H_{i}\left( u,v \right)\mathbb{: = E}\mathrm{\,}\left\lbrack u^{I_{i}\left( a \right)}v^{I_{i}\left( b \right)} \right\rbrack.\]

By the joint distribution above,

\[H_{i}\left( u,v \right) = \frac{\left( n - m_{i} \right)\left( n - m_{i} - 1 \right)}{n\left( n - 1 \right)} + \frac{m_{i}\left( n - m_{i} \right)}{n\left( n - 1 \right)}\left( u + v \right) + \frac{m_{i}\left( m_{i} - 1 \right)}{n\left( n - 1 \right)}\text{uv.}\]

Since

\[L\left( a \right) = \sum_{i = 1}^{T}I_{i}\left( a \right),\ L\left( b \right) = \sum_{i = 1}^{T}I_{i}\left( b \right),\]

we have

\[u^{L\left( a \right)}v^{L\left( b \right)} = \prod_{i = 1}^{T}u^{I_{i}\left( a \right)}v^{I_{i}\left( b \right)}.\]

Because~\(A_{1},\ldots,A_{T}\)~are independently sampled, the random
vectors

\[\left( I_{1}\left( a \right),I_{1}\left( b \right) \right),\ldots,\left( I_{T}\left( a \right),I_{T}\left( b \right) \right)\]

are independent. Hence,

\[\begin{matrix}
 & \begin{matrix}
H\left( u,v \right)\mathbb{\&: = E}\mathrm{\,}\left\lbrack u^{L\left( a \right)}v^{L\left( b \right)} \right\rbrack \\
\& = \prod_{i = 1}^{T}\mathbb{E}\mathrm{\,}\left\lbrack u^{I_{i}\left( a \right)}v^{I_{i}\left( b \right)} \right\rbrack \\
\& = \prod_{i = 1}^{T}H_{i}\left( u,v \right) \\
\& = \prod_{i = 1}^{T}\left\lbrack \frac{\left( n - m_{i} \right)\left( n - m_{i} - 1 \right)}{n\left( n - 1 \right)} + \frac{m_{i}\left( n - m_{i} \right)}{n\left( n - 1 \right)}\left( u + v \right) + \frac{m_{i}\left( m_{i} - 1 \right)}{n\left( n - 1 \right)}\text{uv} \right\rbrack. \\
\end{matrix} & & \mathrm{(10)} \\
\end{matrix}\]

For a bivariate polynomial~\(h\left( u,v \right)\),
let~\(\left\lbrack u^{r}v^{s} \right\rbrack h\left( u,v \right)\)~denote
the coefficient of~\(u^{r}v^{s}\). Then, for
all~\(r,s \in \{ 0,\ldots,T\}\),

\[\mathbb{P}\left( L\left( a \right) = r,L\left( a \right) = s \right) = \left\lbrack u^{r}v^{s} \right\rbrack H\left( u,v \right),\ a \neq b.\]

Recall
\(\mathbf{m} = \left( m_{1},\ldots,m_{T} \right)\text{.\ }\)Define:

\[\begin{matrix}
 & \pi_{r,s}\left( \mathbf{m} \right)\mathbb{: = P(}L\left( a \right) = r,L\left( a \right) = s). & & \mathrm{(11)} \\
\end{matrix}\]

Thus,

\[\begin{matrix}
 & \pi_{r,s}\left( \mathbf{m} \right) = \left\lbrack u^{r}v^{s} \right\rbrack H(u,v). & & \mathrm{(12)} \\
\end{matrix}\]

For~\(S,R \subseteq \left\lbrack T \right\rbrack\), define

\[k_{i}\left( S,R \right): = 1_{\left\{ i \in S \right\}} + 1_{\left\{ i \in R \right\}},\ i \in \left\lbrack T \right\rbrack,\]

and define the two-element MAO pattern weight

\[\begin{matrix}
 & g_{2}(S,R): = \prod_{i = 1}^{T}(\left. \ m_{i}\left. \ ) \right.\ _{k_{i}\left( S,R \right)}\left( n - m_{i} \right)_{2 - k_{i}\left( S,R \right)}. \right.\  & & \mathrm{(13)} \\
\end{matrix}\]

Here, as throughout,

\[\left( x \right)_{q} = x\left( x - 1 \right)\cdots\left( x - q + 1 \right),\ \left( x \right)_{0} = 1,\]

denotes the falling factorial.

The sets~\(S\)~and~\(R\)~respectively record the sampled-set indices in
which~\(a\)~and~\(b\)~occur:

\[S = \left\{ i \in \left\lbrack T \right\rbrack:a \in A_{i} \right\},\ R = \left\{ i \in \left\lbrack T \right\rbrack:b \in A_{i} \right\}.\]

Hence,~\(k_{i}\left( S,R \right) \in \{ 0,1,2\}\)~is the number of
specified elements among~\(a,b\)~that belong to~\(A_{i}\). By
independence across~\(A_{1},\ldots,A_{T}\),

\[\begin{matrix}
 & \mathbb{P}\mathrm{\,}\left( \left\{ i:a \in A_{i} \right\}, = S\mathrm{\,}\left\{ i:b \in A_{i} \right\} = R \right) = \frac{g_{2}\left( S,R \right)}{\left( n \right)_{2}^{T}}. & & \mathrm{(14)} \\
\end{matrix}\]

Consequently,

\[\begin{matrix}
 & \pi_{r,s}\left( \mathbf{m} \right) = \frac{1}{\left( n \right)_{2}^{T}}\sum_{\substack{S,R \subseteq \left\lbrack T \right\rbrack \\ \mid S \mid = r,\mathrm{\ } \mid R \mid = s}}^{}g_{2}\left( S,R \right). & & \mathrm{(15)} \\
\end{matrix}\]

Equations (12) and (15) provide, respectively, the PGF and MAO
pattern-weight representations of the two-element joint occupancy
probability.

\textbf{Theorem 2. Exact cross-second moments}

Let

\[\mu_{t}\mathbb{: = E}\left\lbrack x_{t} \right\rbrack.\]

Then, for all~\(r,s \in \{ 0,\ldots,T\}\),

\[\begin{matrix}
 & \mathbb{E}\left\lbrack x_{r}x_{s} \right\rbrack = 1_{\left\{ r = s \right\}}\mu_{r} + \left( n \right)_{2}\pi_{r,s}\left( \mathbf{m} \right). & & \mathrm{(16)} \\
\end{matrix}\]

Equivalently,

\[\begin{matrix}
 & \mathbb{E}\left\lbrack x_{r}x_{s} \right\rbrack = 1_{\left\{ r = s \right\}}\mu_{r} + \frac{\sum_{\substack{S,R \subseteq \left\lbrack T \right\rbrack \\ \mid S \mid = r,\mathrm{\ } \mid R \mid = s}}^{}g_{2}\left( S,R \right)}{\left( n \right)_{2}^{T - 1}}. & & \mathrm{(17)} \\
\end{matrix}\]

\textbf{Proof.}~By definition,

\[x_{r} = \sum_{a \in \Omega}^{}1_{\left\{ L\left( a \right) = r \right\}},\ x_{s} = \sum_{b \in \Omega}^{}1_{\left\{ L\left( b \right) = s \right\}}.\]

Therefore,

\[x_{r}x_{s} = \left( \sum_{a \in \Omega}^{}1_{\left\{ L\left( a \right) = r \right\}} \right)\left( \sum_{b \in \Omega}^{}1_{\left\{ L\left( b \right) = s \right\}} \right) = \sum_{a \in \Omega}^{}{\sum_{b \in \Omega}^{}{1_{\left\{ L\left( a \right) = r \right\}}1_{\left\{ L\left( b \right) = s \right\}}}}\]

Taking expectations and separating the diagonal and off-diagonal terms
gives

\[E\left( x_{r}x_{s} \right) = \sum_{a \in \Omega}^{}{P(L\left( a \right) = r,L\left( a \right) = s)} + \sum_{\substack{a,b \in \Omega \\ a \neq b}}^{}{P(L\left( a \right) = r,L\left( b \right) = s)}\]

The diagonal probability is zero unless~\(r = s\), in which case it
equals~\(P\left( L\left( a \right) = r \right)\). Hence the diagonal
contribution is

\[1_{\left\{ r = s \right\}}\sum_{a \in \Omega}^{}{P\left( L\left( a \right) = r \right)} = 1_{\left\{ r = s \right\}}\mu_{r}.\]

For~\(a \neq b\), exchangeability of ordered distinct pairs implies that

\[\mathbb{P(}L\left( a \right) = r,L\left( b \right) = s) = \pi_{r,s}\left( \mathbf{m} \right)\]

for every ordered pair~\((a,b)\)~of distinct population elements. Since
there are

\[\left( n \right)_{2} = n\left( n - 1 \right)\]

such ordered pairs, the off-diagonal contribution equals

\[\left( n \right)_{2}\pi_{r,s}\left( \mathbf{m} \right).\]

This proves (16). Equation (17) follows upon substitution of (15).~

\textbf{Corollary 4. Exact covariance structure}

For all~\(r,s \in (0,\ldots,T)\),

\[\begin{matrix}
 & \text{Cov}\left( x_{r},x_{s} \right) = 1_{\left\{ r = s \right\}}\mu_{r} + \left( n \right)_{2}\pi_{r,s}\left( \mathbf{m} \right) - \mu_{r}\mu_{s}. & & \mathrm{(18)} \\
\end{matrix}\]

In particular,

\[\begin{matrix}
 & \text{Var}\left( x_{t} \right) = \mu_{t} + \left( n \right)_{2}\pi_{t,t}\left( \mathbf{m} \right) - \mu_{t}^{2}. & & \mathrm{(19)} \\
\end{matrix}\]

Here,

\[\begin{matrix}
 & \mu_{t} = n\sum_{\substack{S \subseteq \left\lbrack T \right\rbrack \\ \mid S \mid = t}}^{}{\left( \prod_{i \in S}^{}\frac{m_{i}}{n} \right)\left( \prod_{i \notin S}^{}\frac{n - m_{i}}{n} \right)}, & & \mathrm{(20)} \\
\end{matrix}\]

and~\(\pi_{t,t}\left( \mathbf{m} \right)\)~is given by equation (15).
Thus, equation (19) is an exact closed formula for the variance in the
general unequal-size MAO model.

\textbf{Equal-Size Case}

Suppose now that

\[m_{1} = \cdots = m_{T} = m.\]

Define

\[\begin{matrix}
 & \alpha = \frac{m\left( m - 1 \right)}{n\left( n - 1 \right)},\ \beta = \frac{m\left( n - m \right)}{n\left( n - 1 \right)},\ \gamma = \frac{\left( n - m \right)\left( n - m - 1 \right)}{n\left( n - 1 \right)}. & & \mathrm{(21)} \\
\end{matrix}\]

Then

\[\begin{matrix}
 & H\left( u,v \right) = \left\lbrack \gamma + \beta\left( u + v \right) + \alpha uv \right\rbrack^{T}. & & \mathrm{(22)} \\
\end{matrix}\]

Write

\[\pi_{r,s}^{\left( m \right)}\mathbb{: = P(}L\left( a \right) = r,L\left( b \right) = s),\ a \neq b.\]

Let

\[q_{\min} = \max\{ 0,r + s - T\},\ q_{\max} = \min\{ r,s\}.\]

Then

\[\begin{matrix}
 & \pi_{r,s}^{\left( m \right)} = \sum_{q = q_{\min}}^{q_{\max}}\frac{T!}{q!\left( r - q \right)!\left( s - q \right)!\left( T - r - s + q \right)!}\alpha^{q}\beta^{r + s - 2q}\gamma^{T - r - s + q}. & & \mathrm{(23)} \\
\end{matrix}\]

The index~\(q\)~is the number of sampled sets containing
both~\(a\)~and~\(b\). Accordingly, the four membership states occur in

\[q,\ r - q,\ s - q,\ T - r - s + q\]

sampled sets, respectively.

Moreover,

\[\begin{matrix}
 & \mu_{t}^{\left( m \right)}\mathbb{: = E}\left\lbrack x_{t} \right\rbrack = n\left( \frac{T}{t} \right)\left( \frac{m}{n} \right)^{t}\left( 1 - \frac{m}{n} \right)^{T - t}. & & \mathrm{(24)} \\
\end{matrix}\]

Therefore,

\[\begin{matrix}
 & \mathbb{E}\left\lbrack x_{r}x_{s} \right\rbrack = 1_{\left\{ r = s \right\}}\mu_{r}^{\left( m \right)} + n\left( n - 1 \right)\pi_{r,s}^{\left( m \right)}, & & \mathrm{(25)} \\
\end{matrix}\]

\[\begin{matrix}
 & \text{Cov}\left( x_{r},x_{s} \right) = 1_{\left\{ r = s \right\}}\mu_{r}^{\left( m \right)} + n\left( n - 1 \right)\pi_{r,s}^{\left( m \right)} - \mu_{r}^{\left( m \right)}\mu_{s}^{\left( m \right)}, & & \mathrm{(26)} \\
\end{matrix}\]

and, in particular,

\[\begin{matrix}
 & \text{Var}\left( x_{t} \right) = \mu_{t}^{\left( m \right)} + n\left( n - 1 \right)\pi_{t,t}^{\left( m \right)} - \left( \mu_{t}^{\left( m \right)} \right)^{2}. & & \mathrm{(27)} \\
\end{matrix}\]

The joint probability in the variance formula is the single finite sum

\[\begin{matrix}
 & \pi_{t,t}^{\left( m \right)} = \sum_{q = \max\{ 0,2t - T\}}^{t}\frac{T!}{q!\left( t - q \right)!^{2}\left( T - 2t + q \right)!}\alpha^{q}\beta^{2t - 2q}\gamma^{T - 2t + q}. & & \mathrm{(28)} \\
\end{matrix}\]

\textbf{5. Thresholded Occupancy Counts}

For~\(t \in \{ 0,\ldots,T\}\), recall the definition of the thresholded
occupancy count

\[\begin{matrix}
 & x_{\geq t}: = \sum_{q = t}^{T}x_{q} = \sum_{a \in \Omega}^{}1_{\left\{ L\left( a \right) \geq t \right\}}. & & \mathrm{(29)} \\
\end{matrix}\]

Define the one-element upper-tail probability

\[\begin{matrix}
 & p_{\geq t}\left( \mathbf{m} \right)\mathbb{: = P}\left( L\left( a \right) \geq t \right). & & \mathrm{(30)} \\
\end{matrix}\]

By the one-element PGF~\(Q\)

\[\begin{matrix}
 & p_{\geq t}\left( \mathbf{m} \right) = \sum_{q = t}^{T}\lbrack\left. \ z^{q} \right\rbrack Q\left( z \right). & & \mathrm{(31)} \\
\end{matrix}\]

Equivalently,

\[\begin{matrix}
 & p_{\geq t}\left( \mathbf{m} \right) = \sum_{\substack{S \subseteq \left\lbrack T \right\rbrack \\ \mid S \mid \geq t}}^{}{\left( \prod_{i \in S}^{}\frac{m_{i}}{n} \right)\left( \prod_{i \notin S}^{}\frac{n - m_{i}}{n} \right)}. & & \mathrm{(32)} \\
\end{matrix}\]

For two distinct population elements~\(a,b \in \Omega\), define the
joint upper-tail probability

\[\begin{matrix}
 & \Pi_{t,s}\left( \mathbf{m} \right)\mathbb{: = P}\left( L\left( a \right) \geq t,L\left( b \right) \geq s \right). & & \mathrm{(33)} \\
\end{matrix}\]

It follows from the definition
of~\(\pi_{r,s}\left( \mathbf{m} \right)\)~that

\[\begin{matrix}
 & \Pi_{t,s}\left( \mathbf{m} \right) = \sum_{r = t}^{T}{\sum_{q = s}^{T}{\pi_{r,q}\left( \mathbf{m} \right)}}. & & \mathrm{(34)} \\
\end{matrix}\]

Equivalently, by the two-element PGF,

\[\begin{matrix}
 & \Pi_{t,s}\left( \mathbf{m} \right) = \sum_{r = t}^{T}{\sum_{q = s}^{T}{\lbrack\left. \ u^{r}v^{q} \right\rbrack H\left( u,v \right)}}. & & \mathrm{(35)} \\
\end{matrix}\]

Using the two-element MAO pattern weight, we obtain the closed form

\[\begin{matrix}
 & \Pi_{t,s}\left( \mathbf{m} \right) = \frac{1}{\left( n \right)_{2}^{T}}\sum_{\substack{S,R \subseteq \left\lbrack T \right\rbrack \\ \mid S \mid \geq t,\mathrm{\ } \mid R \mid \geq s}}^{}g_{2}\left( S,R \right). & & \mathrm{(36)} \\
\end{matrix}\]

\textbf{Theorem 3. Exact moments of thresholded occupancy counts}

For all~\(t,s \in \{ 0,\ldots,T\}\),

\[\begin{matrix}
 & \mathbb{E}\left\lbrack x_{\geq t} \right\rbrack = n\mathrm{\,}p_{\geq t}\left( \mathbf{m} \right). & & \mathrm{(37)} \\
\end{matrix}\]

Moreover,

\[\begin{matrix}
 & \mathbb{E}\left\lbrack x_{\geq t}x_{\geq s} \right\rbrack = n\mathrm{\,}p_{\geq \max\{ t,s\}}\left( \mathbf{m} \right) + \left( n \right)_{2}\Pi_{t,s}\left( \mathbf{m} \right). & & \mathrm{(38)} \\
\end{matrix}\]

Consequently,

\[\begin{matrix}
 & \text{Cov}\left( x_{\geq t},x_{\geq s} \right) = n\mathrm{\,}p_{\geq \max\{ t,s\}}\left( \mathbf{m} \right) + \left( n \right)_{2}\Pi_{t,s}\left( \mathbf{m} \right) - n^{2}p_{\geq t}\left( \mathbf{m} \right)p_{\geq s}\left( \mathbf{m} \right). & & \mathrm{(39)} \\
\end{matrix}\]

In particular,

\[\begin{matrix}
 & \text{Var}\left( x_{\geq t} \right) = n\mathrm{\,}p_{\geq t}\left( \mathbf{m} \right) + \left( n \right)_{2}\Pi_{t,t}\left( \mathbf{m} \right) - n^{2}p_{\geq t}\left( \mathbf{m} \right)^{2}. & & \mathrm{(40)} \\
\end{matrix}\]

\textbf{Proof.}~Define

\[J_{a}^{\left( \geq t \right)}: = 1_{\left\{ L\left( a \right) \geq t \right\}}.\]

Then

\[x_{\geq t} = \sum_{a \in \Omega}^{}J_{a}^{\left( \geq t \right)}.\]

Therefore,

\[\mathbb{E}\left\lbrack x_{\geq t} \right\rbrack = \sum_{a \in \Omega}^{}\mathbb{P}\left( L\left( a \right) \geq t \right).\]

The MAO model is invariant under permutations of population elements, so
this probability is identical for every~\(a \in \Omega\). Hence,

\[\mathbb{E}\left\lbrack x_{\geq t} \right\rbrack = n\mathrm{\,}p_{\geq t}\left( \mathbf{m} \right),\]

which proves (37).

For the second moment,

\[x_{\geq t}x_{\geq s} = \sum_{a \in \Omega}^{}{\sum_{b \in \Omega}^{}{J_{a}^{\left( \geq t \right)}J_{b}^{\left( \geq s \right)}}}.\]

When~\(a = b\),

\[J_{a}^{\left( \geq t \right)}J_{a}^{\left( \geq s \right)} = 1_{\left\{ L\left( a \right) \geq \max\{ t,s\} \right\}}.\]

Thus, the diagonal contribution is

\[n\mathrm{\,}p_{\geq max\{ t,s\}}\left( \mathbf{m} \right).\]

For~\(a \neq b\), exchangeability of ordered distinct pairs yields

\[\mathbb{P}\left( L\left( a \right) \geq t,L\left( b \right) \geq s \right) = \Pi_{t,s}\left( \mathbf{m} \right).\]

There are~\(\left( n \right)_{2}\)~ordered distinct pairs. Therefore,
the off-diagonal contribution is

\[\left( n \right)_{2}\Pi_{t,s}\left( \mathbf{m} \right).\]

This proves (38). Equations (39) and (40) follow by subtracting the
appropriate products of means.

\textbf{Equal-Size Case}

Suppose that

\[m_{1} = \cdots = m_{T} = m.\]

Define

\[\begin{matrix}
 & p_{\geq t}^{\left( m \right)}\mathbb{: = P}\left( L\left( a \right) \geq t \right) = \sum_{r = t}^{T}{\left( \frac{T}{r} \right)\left( \frac{m}{n} \right)^{r}\left( 1 - \frac{m}{n} \right)^{T - r}}. & & \mathrm{(41)} \\
\end{matrix}\]

and

\[\begin{matrix}
 & \Pi_{t,s}^{\left( m \right)}\mathbb{: = P}\left( L\left( a \right) \geq t,L\left( b \right) \geq s \right). & & \mathrm{(42)} \\
\end{matrix}\]

Then

\[\begin{matrix}
 & \Pi_{t,s}^{\left( m \right)} = \sum_{r = t}^{T}{\sum_{q = s}^{T}\pi_{r,q}^{\left( m \right)}}, & & \mathrm{(43)} \\
\end{matrix}\]

where~\(\pi_{r,q}^{\left( m \right)}\)~is given by equation (23).

Consequently,

\[\begin{matrix}
 & \mathbb{E}\left\lbrack x_{\geq t} \right\rbrack = n\mathrm{\,}p_{\geq t}^{\left( m \right)}, & & \mathrm{(44)} \\
\end{matrix}\]

\[\begin{matrix}
 & \mathbb{E}\left\lbrack x_{\geq t}x_{\geq s} \right\rbrack = n\mathrm{\,}p_{\geq \max\{ t,s\}}^{\left( m \right)} + n\left( n - 1 \right)\Pi_{t,s}^{\left( m \right)}, & & \mathrm{(45)} \\
\end{matrix}\]

\[\begin{matrix}
 & \text{Cov}\left( x_{\geq t},x_{\geq s} \right) = n\mathrm{\,}p_{\geq max\{ t,s\}}^{\left( m \right)} + n\left( n - 1 \right)\Pi_{t,s}^{\left( m \right)} - n^{2}p_{\geq t}^{\left( m \right)}p_{\geq s}^{\left( m \right)}, & & \mathrm{(46)} \\
\end{matrix}\]

and

\[\begin{matrix}
 & \text{Var}\left( x_{\geq t} \right) = n\mathrm{\,}p_{\geq t}^{\left( m \right)} + n\left( n - 1 \right)\Pi_{t,t}^{\left( m \right)} - n^{2}\left( p_{\geq t}^{\left( m \right)} \right)^{2}. & & \mathrm{(47)} \\
\end{matrix}\]

\textbf{6. Higher-Order Joint PGFs and MAO Factorial Moments}

The one-element PGFs~\(Q_{i}\)~and~\(Q\), and the two-element
PGFs~\(H_{i}\)~and~\(H\), are the first two cases of a natural hierarchy
of joint PGFs. To distinguish the dimension of a joint PGF from the set
index, we use the superscript~\(\left( \mathcal{l} \right)\)~for the
number of population elements and retain the subscript~\(i\)~for the
sampled set~\(A_{i}\). Thus,

\[Q_{i} = \Phi_{i}^{\left( 1 \right)},\ Q = \Phi^{\left( 1 \right)},\ H_{i} = \Phi_{i}^{\left( 2 \right)},\ H = \Phi^{\left( 2 \right)}.\]

Fix an
ordered~\(\mathcal{l}\)-tuple~\(\left( a_{1},\ldots,a_{\mathcal{l}} \right)\)~of
pairwise distinct elements of~\(\Omega\),
where~\(1 \leq \mathcal{l} \leq n\). For
each~\(i \in \left\lbrack T \right\rbrack\), define the
per-set~\(\mathcal{l}\)-element joint PGF by

\[\Phi_{i}^{\left( \mathcal{l} \right)}\left( u_{1},\ldots,u_{\mathcal{l}} \right)\mathbb{: = E}\mathrm{\,}\left\lbrack \prod_{j = 1}^{\mathcal{l}}u_{j}^{I_{i}\left( a_{j} \right)} \right\rbrack.\]

Define the corresponding aggregate~\(\mathcal{l}\)-element joint PGF by

\[\begin{matrix}
 & \Phi^{\left( \mathcal{l} \right)}\left( u_{1},\ldots,u_{\mathcal{l}} \right)\mathbb{: = E}\mathrm{\,}\left\lbrack \prod_{j = 1}^{\mathcal{l}}u_{j}^{L\left( a_{j} \right)} \right\rbrack. & & \mathrm{(48)} \\
\end{matrix}\]

For~\(q_{1},\ldots,q_{\mathcal{l}} \in \{ 0,\ldots,T\}\),

\[\begin{matrix}
 & \left\lbrack u_{1}^{q_{1}}\cdots u_{\mathcal{l}}^{q_{\mathcal{l}}} \right\rbrack\Phi^{\left( \mathcal{l} \right)}\left( u_{1},\ldots,u_{\mathcal{l}} \right)\mathbb{= P}\left( L\left( a_{1} \right) = q_{1},\ldots,L\left( a_{\mathcal{l}} \right) = q_{\mathcal{l}} \right). & & \mathrm{(49)} \\
\end{matrix}\]

\textbf{Theorem 4.~}\(\mathcal{l}\)\textbf{-element joint probability
generating function}

For~\(1 \leq \mathcal{l} \leq n\),

\[\begin{matrix}
 & \begin{matrix}
\Phi^{\left( \mathcal{l} \right)}\left( u_{1},\ldots,u_{\mathcal{l}} \right)\& = \prod_{i = 1}^{T}\Phi_{i}^{\left( \mathcal{l} \right)}\left( u_{1},\ldots,u_{\mathcal{l}} \right) \\
\& = \prod_{i = 1}^{T}\left\lbrack \sum_{B \subseteq \left\lbrack \mathcal{l} \right\rbrack}^{}\frac{\left( m_{i} \right)_{\mid B \mid}\left( n - m_{i} \right)_{\mathcal{l} - \mid B \mid}}{\left( n \right)_{\mathcal{l}}}\prod_{j \in B}^{}u_{j} \right\rbrack. \\
\end{matrix} & & \mathrm{(50)} \\
\end{matrix}\]

Here,

\[\left( x \right)_{k} = x\left( x - 1 \right)\cdots\left( x - k + 1 \right),\ \left( x \right)_{0} = 1,\]

denotes the falling factorial.

\textbf{Proof.}~Since

\[L\left( a_{j} \right) = \sum_{i = 1}^{T}I_{i}\left( a_{j} \right),\]

we have

\[\prod_{j = 1}^{\mathcal{l}}u_{j}^{L\left( a_{j} \right)} = \prod_{i = 1}^{T}{\prod_{j = 1}^{\mathcal{l}}u_{j}^{I_{i}\left( a_{j} \right)}}.\]

The random vectors

\[\left( I_{i}\left( a_{1} \right),\ldots,I_{i}\left( a_{\mathcal{l}} \right) \right),\ i \in \left\lbrack T \right\rbrack,\]

are independent because~\(A_{1},\ldots,A_{T}\)~are independent.
Therefore,

\[\Phi^{\left( \mathcal{l} \right)}\left( u_{1},\ldots,u_{\mathcal{l}} \right) = \prod_{i = 1}^{T}\mathbb{E}\mathrm{\,}\left\lbrack \prod_{j = 1}^{\mathcal{l}}u_{j}^{I_{i}\left( a_{j} \right)} \right\rbrack = \prod_{i = 1}^{T}\Phi_{i}^{\left( \mathcal{l} \right)}\left( u_{1},\ldots,u_{\mathcal{l}} \right).\]

Fix~\(i \in \left\lbrack T \right\rbrack\), and let

\[B = \left\{ j \in \left\lbrack \mathcal{l} \right\rbrack:a_{j} \in A_{i} \right\}.\]

If~\(\mid B \mid = k\), then precisely the~\(k\)~specified elements
indexed by~\(B\)~are included in~\(A_{i}\), while the
remaining~\(\mathcal{l} - k\)~specified elements are excluded. The
probability of this membership pattern is

\[\frac{\left( m_{i} \right)_{k}\left( n - m_{i} \right)_{\mathcal{l} - k}}{\left( n \right)_{\mathcal{l}}},\]

and its monomial contribution is

\[\prod_{j \in B}^{}u_{j}.\]

Consequently,

\[\Phi_{i}^{\left( \mathcal{l} \right)}\left( u_{1},\ldots,u_{\mathcal{l}} \right) = \sum_{B \subseteq \left\lbrack \mathcal{l} \right\rbrack}^{}{\frac{\left( m_{i} \right)_{\mid B \mid}\left( n - m_{i} \right)_{\mathcal{l} - \mid B \mid}}{\left( n \right)_{\mathcal{l}}}\prod_{j \in B}^{}u_{j}}.\]

Summing over
all~\(B \subseteq \left\lbrack \mathcal{l} \right\rbrack\)~gives the
per-set factor~\(\Phi_{i}^{\left( \mathcal{l} \right)}\), and
multiplication over~\(i = 1,\ldots,T\)~proves (50).

\textbf{Definition 2. MAO pattern weight}

Let
\(S_{1},\ldots,S_{\mathcal{l}} \subseteq \left\lbrack T \right\rbrack\)
be the membership patterns of~\(\mathcal{l}\)~distinct population
elements. For~\(i \in \left\lbrack T \right\rbrack\), define

\[\begin{matrix}
 & k_{i}\left( S_{1},\ldots,S_{\mathcal{l}} \right): = \sum_{j = 1}^{\mathcal{l}}1_{\left\{ i \in S_{j} \right\}}. & & \mathrm{(51)} \\
\end{matrix}\]

Define the~\(\mathcal{l}\)-element MAO pattern weight by

\[\begin{matrix}
 & g_{\mathcal{l}}\left( S_{1},\ldots,S_{\mathcal{l}} \right): = \prod_{i = 1}^{T}(\left. \ m_{i}\left. \ ) \right.\ _{k_{i}\left( S_{1},\ldots,S_{\mathcal{l}} \right)}\left( n - m_{i} \right)_{\mathcal{l} - k_{i}\left( S_{1},\ldots,S_{\mathcal{l}} \right)}. \right.\  & & \mathrm{(52)} \\
\end{matrix}\]

Thus,~\(k_{i}\left( S_{1},\ldots,S_{\mathcal{l}} \right)\)~is the number
of the~\(\mathcal{l}\)~specified population elements that belong
to~\(A_{i}\)~under the pattern
tuple~\(\left( S_{1},\ldots,S_{\mathcal{l}} \right)\).

\textbf{Proposition 2. Probabilistic interpretation of the MAO pattern
weight}

For pairwise distinct~\(a_{1},\ldots,a_{\mathcal{l}} \in \Omega\),

\[\begin{matrix}
 & \mathbb{P}\left( \left\{ i \in \left\lbrack T \right\rbrack:a_{j} \in A_{i} \right\} = S_{j}\mathrm{\text{~for~all~}}j \in \left\lbrack \mathcal{l} \right\rbrack \right) = \frac{g_{\mathcal{l}}\left( S_{1},\ldots,S_{\mathcal{l}} \right)}{\left( n \right)_{\mathcal{l}}^{T}}. & & \mathrm{(53)} \\
\end{matrix}\]

\textbf{Proof.}~For each~\(i \in \left\lbrack T \right\rbrack\), exactly

\[k_{i}\left( S_{1},\ldots,S_{\mathcal{l}} \right)\]

of the~\(\mathcal{l}\)~specified elements must belong to~\(A_{i}\),
while the remaining

\[\mathcal{l} - k_{i}\left( S_{1},\ldots,S_{\mathcal{l}} \right)\]

specified elements must not belong to~\(A_{i}\). The probability of this
membership pattern for~\(A_{i}\)~is

\[\frac{\left( m_{i} \right)_{k_{i}\left( S_{1},\ldots,S_{\mathcal{l}} \right)}\left( n - m_{i} \right)_{\mathcal{l} - k_{i}\left( S_{1},\ldots,S_{\mathcal{l}} \right)}}{\left( n \right)_{\mathcal{l}}}.\]

Independence of~\(A_{1},\ldots,A_{T}\)~gives equation (53).

For~\(\mathcal{l} = 1\),

\[g_{1}\left( S \right) = \prod_{i \in S}^{}m_{i}\prod_{i \notin S}^{}(\left. \ n - m_{i} \right),\]

which is the original one-element MAO weight.

\textbf{Definition 3. Generalized MAO factorial moments}

For~\(1 \leq \mathcal{l} \leq n\), let

\[B_{1},\ldots,B_{\mathcal{l}} \subseteq \{ 0,1,\ldots,T\}\]

be nonempty sets of admissible occupancy levels. Define the
generalized~\(\mathcal{l}\)-th MAO factorial moment by

\[\begin{matrix}
 & F_{\mathcal{l}}\left( B_{1},\ldots,B_{\mathcal{l}} \right)\mathbb{: = E}\mathrm{\,}\left\lbrack \sum_{\substack{\left( a_{1},\ldots,a_{\mathcal{l}} \right) \in \Omega^{\mathcal{l}} \\ a_{1},\ldots,a_{\mathcal{l}}\mathrm{\ }\text{pairwise}\mathrm{\ }\text{distinct}}}^{}{\prod_{j = 1}^{\mathcal{l}}1_{\left\{ L\left( a_{j} \right) \in B_{j} \right\}}} \right\rbrack. & & \mathrm{(54)} \\
\end{matrix}\]

Equivalently,

\[\begin{matrix}
 & F_{\mathcal{l}}\left( B_{1},\ldots,B_{\mathcal{l}} \right) = \sum_{\substack{\left( a_{1},\ldots,a_{\mathcal{l}} \right) \in \Omega^{\mathcal{l}} \\ a_{1},\ldots,a_{\mathcal{l}}\mathrm{\ }\text{pairwise}\mathrm{\ }\text{distinct}}}^{}\mathbb{P}\mathrm{\,}\left( L\left( a_{j} \right) \in B_{j}\mathrm{\text{~for~all~}}j \in \left\lbrack \mathcal{l} \right\rbrack \right). & & \mathrm{(55)} \\
\end{matrix}\]

By exchangeability of the population labels, for any fixed pairwise
distinct elements~\(a_{1},\ldots,a_{\mathcal{l}} \in \Omega\),

\[\begin{matrix}
 & F_{\mathcal{l}}\left( B_{1},\ldots,B_{\mathcal{l}} \right) = \left( n \right)_{\mathcal{l}}\mathrm{\,}\mathbb{P}\mathrm{\,}\left( L\left( a_{j} \right) \in B_{j}\mathrm{\text{~for~all~}}j \in \left\lbrack \mathcal{l} \right\rbrack \right). & & \mathrm{(56)} \\
\end{matrix}\]

Equivalently, in terms of the~\(\mathcal{l}\)-element joint probability
generating function,

\[\begin{matrix}
 & F_{\mathcal{l}}\left( B_{1},\ldots,B_{\mathcal{l}} \right) = \left( n \right)_{\mathcal{l}}\sum_{q_{1} \in B_{1}}^{}\cdots\sum_{q_{\mathcal{l}} \in B_{\mathcal{l}}}^{}\lbrack\left. \ u_{1}^{q_{1}}\cdots u_{\mathcal{l}}^{q_{\mathcal{l}}} \right\rbrack\Phi^{\left( \mathcal{l} \right)}\left( u_{1},\ldots,u_{\mathcal{l}} \right). & & \mathrm{(57)} \\
\end{matrix}\]

The terminology reflects the fact that, when all admissible-level sets
are identical:

\[B_{1} = \cdots = B_{\mathcal{l}} = B,\]

we have

\[\begin{matrix}
 & F_{\mathcal{l}}\left( B,\ldots,B \right)\mathbb{= E}\left\lbrack \left( \sum_{a \in \Omega}^{}1_{\left\{ L\left( a \right) \in B \right\}} \right)_{\mathcal{l}} \right\rbrack\mathrm{\,}. & & \mathrm{(58)} \\
\end{matrix}\]

Thus,~\(F_{\mathcal{l}}\left( B,\ldots,B \right)\)~is
the~\(\mathcal{l}\)-th factorial moment of the number of population
elements whose occupancy levels belong to~\(B\).

For an exact occupancy level~\(t \in \{ 0,1,\ldots,T\}\), define the
diagonal specialization

\[\begin{matrix}
 & F_{\mathcal{l}}\left( t \right): = F_{\mathcal{l}}\left( \left\{ t \right\},\ldots,\{ t\} \right)\mathbb{= E}\left\lbrack \left( x_{t} \right)_{\mathcal{l}} \right\rbrack. & & \mathrm{(59)} \\
\end{matrix}\]

For the thresholded count, let

\[B_{\geq t}: = \{ t,t + 1,\ldots,T\},\]

define

\[\begin{matrix}
 & F_{\mathcal{l}}^{\geq}\left( t \right): = F_{\mathcal{l}}\left( B_{\geq t},\ldots,B_{\geq t} \right)\mathbb{= E}\left\lbrack \left( x_{\geq t} \right)_{\mathcal{l}} \right\rbrack. & & \mathrm{(60)} \\
\end{matrix}\]

Thus, the exact and thresholded MAO factorial moments are diagonal
specializations of the generalized MAO factorial-moment functional.

\textbf{Definition 4. Generalized MAO transversal sum}

For nonempty sets

\[B_{1},\ldots,B_{\mathcal{l}} \subseteq \{ 0,1,\ldots,T\},\]

define

\[\begin{matrix}
 & G_{T}\left( B_{1},\ldots,B_{\mathcal{l}} \right): = \sum_{\substack{S_{1},\ldots,S_{\mathcal{l}} \subseteq \left\lbrack T \right\rbrack \\ \mid S_{j} \mid \in B_{j},\mathrm{\ }j \in \left\lbrack \mathcal{l} \right\rbrack}}^{}g_{\mathcal{l}}\left( S_{1},\ldots,S_{\mathcal{l}} \right). & & \mathrm{(61)} \\
\end{matrix}\]

For singleton admissible-level sets, use the abbreviation

\[\begin{matrix}
 & G_{T}\left( p_{1},\ldots,p_{\mathcal{l}} \right): = G_{T}\left( \left\{ p_{1} \right\},\ldots,\left\{ p_{\mathcal{l}} \right\} \right). & & \mathrm{(62)} \\
\end{matrix}\]

Partitioning the transversal sum according to the pattern cardinalities
gives

\[\begin{matrix}
 & G_{T}\left( B_{1},\ldots,B_{\mathcal{l}} \right) = \sum_{p_{1} \in B_{1}}^{}\cdots\sum_{p_{\mathcal{l}} \in B_{\mathcal{l}}}^{}G_{T}\left( p_{1},\ldots,p_{\mathcal{l}} \right). & & \mathrm{(63)} \\
\end{matrix}\]

Specifically, for a nonempty set~\(B \subseteq \{ 0,1,\ldots,T\}\),
write

\[B^{\mathcal{l}}: = \underset{\mathcal{l}\mathrm{\text{~times}}}{}\]

for the diagonal~\(\mathcal{l}\)-tuple or \(p \in \{ 0,1,\ldots,T\}\),
abbreviate

\[G_{T}\left( B^{\mathcal{l}} \right): = G_{T}\left( B,\ldots,B \right),\ F_{\mathcal{l}}\left( B^{\mathcal{l}} \right): = F_{\mathcal{l}}\left( B,\ldots,B \right)\]

\[G_{T}\left( p^{\mathcal{l}} \right): = G_{T}\left( \left\{ p \right\},\ \ldots,\{ p\} \right),\ F_{\mathcal{l}}\left( p^{\mathcal{l}} \right): = F_{\mathcal{l}}\left( \left\{ p \right\},\ \ldots,\{ p\} \right)\]

where each right-hand side has~\(\mathcal{l}\)~copies of~\(B\) or
\(\{ p\}\). Then the above defined
\(F_{\mathcal{l}}\left( t \right) = F_{\mathcal{l}}\left( {\{ t\}}^{\mathcal{l}} \right),\ F_{\mathcal{l}}^{\geq}\left( t \right) = F_{\mathcal{l}}\left( B_{\geq t}^{\mathcal{l}} \right)\).

\textbf{Theorem 5. Generalized exact MAO factorial-moment
representation}

For every~\(1 \leq \mathcal{l} \leq n\)~and all nonempty sets

\[B_{1},\ldots,B_{\mathcal{l}} \subseteq \{ 0,1,\ldots,T\},\]

the generalized MAO factorial-moment functional has the exact
representation

\[\begin{matrix}
 & F_{\mathcal{l}}\left( B_{1},\ldots,B_{\mathcal{l}} \right) = \frac{G_{T}\left( B_{1},\ldots,B_{\mathcal{l}} \right)}{\left( n \right)_{\mathcal{l}}^{T - 1}}. & & \mathrm{(64)} \\
\end{matrix}\]

Equivalently,

\[\begin{matrix}
 & F_{\mathcal{l}}\left( B_{1},\ldots,B_{\mathcal{l}} \right) = \frac{1}{\left( n \right)_{\mathcal{l}}^{T - 1}}\sum_{\substack{S_{1},\ldots,S_{\mathcal{l}} \subseteq \left\lbrack T \right\rbrack \\ \mid S_{j} \mid \in B_{j},\mathrm{\ }j \in \left\lbrack \mathcal{l} \right\rbrack}}^{}{\prod_{i = 1}^{T}{\left( m_{i} \right)_{k_{i}\left( S_{1},..,S_{\mathcal{l}} \right)}\left( n - m_{i} \right)_{\mathcal{l} - k_{i}\left( S_{1},..,S_{\mathcal{l}} \right)}}} & & \mathrm{(65)} \\
\end{matrix}\]

The corresponding joint-PGF representation is

\[\begin{matrix}
 & F_{\mathcal{l}}\left( B_{1},\ldots,B_{\mathcal{l}} \right) = \left( n \right)_{\mathcal{l}}\sum_{q_{1} \in B_{1}}^{}\cdots\sum_{q_{\mathcal{l}} \in B_{\mathcal{l}}}^{}{\left\lbrack u_{1}^{q_{1}}\cdots u_{\mathcal{l}}^{q_{\mathcal{l}}} \right\rbrack\Phi^{\left( \mathcal{l} \right)}\left( u_{1},\ldots,u_{\mathcal{l}} \right)}. & & \mathrm{(66)} \\
\end{matrix}\]

Moreover,

\[\begin{matrix}
 & F_{\mathcal{l}}\left( B_{1},\ldots,B_{\mathcal{l}} \right) = \sum_{p_{1} \in B_{1}}^{}\cdots\sum_{p_{\mathcal{l}} \in B_{\mathcal{l}}}^{}{F_{\mathcal{l}}\left( p_{1},\ldots,p_{\mathcal{l}} \right)}. & & \mathrm{(67)} \\
\end{matrix}\]

\textbf{Proof.}~Fix pairwise distinct
elements~\(a_{1},\ldots,a_{\mathcal{l}} \in \Omega\), and define their
membership patterns by

\[S_{j}: = \left\{ i \in \left\lbrack T \right\rbrack:a_{j} \in A_{i} \right\},\ j \in \left\lbrack \mathcal{l} \right\rbrack.\]

Since~\(L\left( a_{j} \right) = \mid S_{j} \mid\), the event

\[\left\{ L\left( a_{j} \right) \in B_{j}\mathrm{\text{~for~every~}}j \in \left\lbrack \mathcal{l} \right\rbrack \right\}\]

is the disjoint union of the membership-pattern events indexed by
tuples~\(\left( S_{1},..,S_{\mathcal{l}} \right)\)~satisfying

\[\mid S_{j} \mid \in B_{j},\ j \in \left\lbrack \mathcal{l} \right\rbrack.\]

Proposition 2 therefore gives

\[\mathbb{P}\left( L\left( a_{j} \right) \in B_{j}\mathrm{\text{~for~every~}}j \in \left\lbrack \mathcal{l} \right\rbrack \right) = \frac{G_{T}\left( B_{1},\ldots,B_{\mathcal{l}} \right)}{\left( n \right)_{\mathcal{l}}^{T}}.\]

Multiplication by~\(\left( n \right)_{\mathcal{l}}\), as in (56), proves
(64) and (65). Equation (66) follows from the coefficient-extraction
identity (49), and (67) follows by partitioning the transversal sum
according
to~\(\left( \left\lfloor S_{1} \right\rfloor,\ldots,\left\lfloor S_{\mathcal{l}} \right\rfloor \right)\).

\textbf{Remark 1.}~The summation in (65) permits coincident membership
patterns; that is,~\(S_{j} = S_{k}\)~may hold for~\(j \neq k\). The
population elements~\(a_{1},\ldots,a_{\mathcal{l}}\), however, remain
pairwise distinct, as required by the falling factorial. This
distinguishes the factorial-moment expansion from an ordinary-power
expansion, in which repeated population elements also occur.

\textbf{Remark 2.}~The generalized MAO norm introduced in (Mao and Xue,
2025) is precisely a pointwise MAO factorial moment:

\[\begin{matrix}
 & N_{T}\left( p_{1},\ldots,p_{\mathcal{l}} \right) = F_{\mathcal{l}}\left( p_{1},\ldots,p_{\mathcal{l}} \right). & & \mathrm{(68)} \\
\end{matrix}\]

The notation~\(F_{\mathcal{l}}\)~is used in the present paper to make
this probabilistic interpretation explicit.

\textbf{Corollary 5. Exact-occupancy factorial moments}

For every~\(t \in \{ 0,1,\ldots,T\}\)~and~\(1 \leq \mathcal{l} \leq n\),

\[\begin{matrix}
 & F_{\mathcal{l}}\left( t \right)\mathbb{= E}\left\lbrack \left( x_{t} \right)_{\mathcal{l}} \right\rbrack = \frac{G_{T}\left( t,\ldots,t \right)}{\left( n \right)_{\mathcal{l}}^{T - 1}} = \frac{G_{T}\left( t^{\mathcal{l}} \right)}{\left( n \right)_{\mathcal{l}}^{T - 1}}. & & \mathrm{(69)} \\
\end{matrix}\]

Equivalently,

\[\begin{matrix}
 & F_{\mathcal{l}}\left( t \right) = \frac{1}{\left( n \right)_{\mathcal{l}}^{T - 1}}\sum_{\substack{S_{1},\ldots,S_{\mathcal{l}} \subseteq \left\lbrack T \right\rbrack \\ \mid S_{1} \mid = \cdots = \mid S_{\mathcal{l}} \mid = t}}^{}{g_{\mathcal{l}}\left( S_{1},..,S_{\mathcal{l}} \right)}. & & \mathrm{(70)} \\
\end{matrix}\]

\textbf{Proof.}~Apply Theorem 5 with
\(B_{1} = \cdots = B_{\mathcal{l}} = \left\{ t \right\}.\)

\textbf{Corollary 6. Threshold-occupancy factorial moments}

For every~\(t \in \{ 0,1,\ldots,T\}\)~and~\(1 \leq \mathcal{l} \leq n\),

\[\begin{matrix}
 & F_{\mathcal{l}}^{\geq}\left( t \right)\mathbb{= E}\left\lbrack \left( x_{\geq t} \right)_{\mathcal{l}} \right\rbrack = \frac{G_{T}\left( B_{\geq t},\ldots,B_{\geq t} \right)}{\left( n \right)_{\mathcal{l}}^{T - 1}} = \frac{G_{T}\left( B_{\geq t}^{\mathcal{l}} \right)}{\left( n \right)_{\mathcal{l}}^{T - 1}}. & & \mathrm{(71)} \\
\end{matrix}\]

Equivalently,

\[\begin{matrix}
 & F_{\mathcal{l}}^{\geq}\left( t \right) = \frac{1}{\left( n \right)_{\mathcal{l}}^{T - 1}}\sum_{\substack{S_{1},\ldots,S_{\mathcal{l}} \subseteq \left\lbrack T \right\rbrack \\ \mid S_{j} \mid \geq t,\mathrm{\ }j \in \left\lbrack \mathcal{l} \right\rbrack}}^{}g_{\mathcal{l}}\left( S_{1},..,S_{\mathcal{l}} \right). & & \mathrm{(72)} \\
\end{matrix}\]

Equivalently, in terms of the pointwise MAO factorial moments,

\[\begin{matrix}
 & F_{\mathcal{l}}^{\geq}\left( t \right) = \sum_{q_{1} = t}^{T}\cdots\sum_{q_{\mathcal{l}} = t}^{T}F_{\mathcal{l}}\left( q_{1},\ldots,q_{\mathcal{l}} \right). & & \mathrm{(73)} \\
\end{matrix}\]

\textbf{Proof.}~Apply Theorem 5 and (67) with
\(B_{1} = \cdots = B_{\mathcal{l}} = B_{\geq t} = \{ t,t + 1,\ldots,T\}.\)

\textbf{Corollary 7. Exact raw moments}

For every~\(t \in \{ 0,1,\ldots,T\}\)~and every positive integer~\(v\),

\[\begin{matrix}
 & \mathbb{E}\left\lbrack x_{t}^{v} \right\rbrack = \sum_{\mathcal{l} = 1}^{\min\{ v,n\}}\left\{ \frac{v}{\mathcal{l}} \right\} F_{\mathcal{l}}\left( t \right), & & \mathrm{(74)} \\
\end{matrix}\]

and

\[\begin{matrix}
 & \mathbb{E}\left\lbrack x_{\geq t}^{v} \right\rbrack = \sum_{\mathcal{l} = 1}^{\min\{ v,n\}}\left\{ \frac{v}{\mathcal{l}} \right\} F_{\mathcal{l}}^{\geq}\left( t \right), & & \mathrm{(75)} \\
\end{matrix}\]

where

\[\left\{ \frac{v}{\mathcal{l}} \right\}\]

denotes a Stirling number of the second kind. Consequently,

\[\begin{matrix}
 & \mathbb{E}\left\lbrack x_{t}^{v} \right\rbrack = \sum_{\mathcal{l} = 1}^{\min\{ v,n\}}\left\{ \frac{v}{\mathcal{l}} \right\}\frac{G_{T}\left( t^{\mathcal{l}} \right)}{\left( n \right)_{\mathcal{l}}^{T - 1}}, & & \mathrm{(76)} \\
\end{matrix}\]

and

\[\begin{matrix}
 & \mathbb{E}\left\lbrack x_{\geq t}^{v} \right\rbrack = \sum_{\mathcal{l} = 1}^{\min\{ v,n\}}\left\{ \frac{v}{\mathcal{l}} \right\}\frac{G_{T}\left( B_{\geq t}^{\mathcal{l}} \right)}{\left( n \right)_{\mathcal{l}}^{T - 1}}. & & \mathrm{(77)} \\
\end{matrix}\]

An equivalent pointwise representation of the threshold raw moment is

\[\begin{matrix}
 & \mathbb{E}\left\lbrack x_{\geq t}^{v} \right\rbrack = \sum_{\mathcal{l} = 1}^{\min\{ v,n\}}\left\{ \frac{v}{\mathcal{l}} \right\}\sum_{q_{1} = t}^{T}\cdots\sum_{q_{\mathcal{l}} = t}^{T}{F_{\mathcal{l}}\left( q_{1},\ldots,q_{\mathcal{l}} \right)}. & & \mathrm{(78)} \\
\end{matrix}\]

\textbf{Proof.}~For every nonnegative integer~\(x\),

\[\begin{matrix}
 & x^{v} = \sum_{\mathcal{l} = 0}^{v}\left\{ \frac{v}{\mathcal{l}} \right\}\left( x \right)_{\mathcal{l}}. & & \mathrm{(79)} \\
\end{matrix}\]

Applying (79) to~\(x_{t}\)~and~\(x_{\geq t}\), and then taking
expectations, gives (74) and (75). Since both counts are bounded above
by~\(n\), all terms with~\(\mathcal{l} > n\)~vanish. Since~\(v \geq 1\),
the term with~\(\mathcal{l} = 0\)~also vanishes. Substitution of
Corollaries 5 and 6 gives (76)--(78).~

Therefore, exact- and threshold-occupancy moments are obtained from the
same generalized MAO factorial-moment functional by diagonal
specialization and Stirling transformation.

\textbf{7. Fully Occupied Counts and Factorial-Moment Inequalities}

The fully occupied count is

\[x_{T} = \left. \mid\bigcap_{i = 1}^{T}A_{i} \right.\mid.\]

Its mean is

\[\begin{matrix}
 & \mu_{T}\mathbb{: = E}\left\lbrack x_{T} \right\rbrack = \frac{\prod_{i = 1}^{T}m_{i}}{n^{T - 1}}. & & \mathrm{(80)} \\
\end{matrix}\]

\textbf{Theorem 6. Fully occupied MAO factorial-moment inequality}

For every integer~\(1 \leq \mathcal{l} \leq n\),

\[\begin{matrix}
 & F_{\mathcal{l}}\left( T \right)\mathbb{= E}\left\lbrack \left( x_{T} \right)_{\mathcal{l}} \right\rbrack = \frac{\prod_{i = 1}^{T}(\left. \ m_{i}\left. \ ) \right.\ _{\mathcal{l}} \right.\ }{\left( n \right)_{\mathcal{l}}^{T - 1}}. & & \mathrm{(81)} \\
\end{matrix}\]

Moreover,

\[\begin{matrix}
 & F_{\mathcal{l}}\left( T \right) \leq \mu_{T}^{\mathcal{l}}. & & \mathrm{(82)} \\
\end{matrix}\]

If~\(\mu_{T} > 0\)~and~\(2 \leq \mathcal{l} \leq n\), then

\[\begin{matrix}
 & F_{\mathcal{l}}\left( T \right) < \mu_{T}^{\mathcal{l}}. & & \mathrm{(83)} \\
\end{matrix}\]

Thus, every factorial moment of order at least two is strictly bounded
above by the corresponding power of the mean whenever~\(\mu_{T} > 0\).

\textbf{Proof.}~For the specialization

\[B_{1} = \cdots = B_{\mathcal{l}} = \left\{ T \right\},\]

the only admissible membership-pattern tuple in Theorem 5 is

\[S_{1} = \cdots = S_{\mathcal{l}} = \left\lbrack T \right\rbrack.\]

Therefore,

\[k_{i}\left( S_{1},..,S_{\mathcal{l}} \right) = \mathcal{l}\ \mathrm{\text{for~every~}}i \in \left\lbrack T \right\rbrack,\]

and hence

\[G_{T}\left( T,\ldots,T \right) = \prod_{i = 1}^{T}(\left. \ m_{i}\left. \ ) \right.\ _{\mathcal{l}}. \right.\ \]

The generalized MAO representation in Theorem 5 now gives (81).

For an integer~\(s \geq \mathcal{l}\), define

\[\begin{matrix}
 & r_{\mathcal{l}}\left( s \right): = \frac{\left( s \right)_{\mathcal{l}}}{s^{\mathcal{l}}} = \prod_{j = 0}^{\mathcal{l} - 1}\left( 1 - \frac{j}{s} \right). & & \mathrm{(84)} \\
\end{matrix}\]

For fixed~\(\mathcal{l}\), the
function~\(r_{\mathcal{l}}\left( s \right)\)~is nondecreasing for
integers~\(s \geq \mathcal{l}\). If some~\(m_{i} < \mathcal{l}\), then

\[F_{\mathcal{l}}\left( T \right) = 0,\]

and (82) follows immediately. Otherwise,~\(m_{i} \geq \mathcal{l}\)~for
every~\(i\), and~\(m_{i} \leq n\)~implies

\[r_{\mathcal{l}}\left( m_{i} \right) \leq r_{\mathcal{l}}\left( n \right).\]

Using (80) and (81), we obtain

\[\begin{matrix}
 & \frac{F_{\mathcal{l}}\left( T \right)}{\mu_{T}^{\mathcal{l}}} = \frac{\prod_{i = 1}^{T}r_{\mathcal{l}}\left( m_{i} \right)}{r_{\mathcal{l}}\left( n \right)^{T - 1}} \leq r_{\mathcal{l}}\left( n \right) \leq 1. & & \mathrm{(85)} \\
\end{matrix}\]

This proves (82). If~\(\mu_{T} > 0\)~and~\(2 \leq \mathcal{l} \leq n\),
then either some~\(m_{i} < \mathcal{l}\), in which case

\[F_{\mathcal{l}}\left( T \right) = 0 < \mu_{T}^{\mathcal{l}},\]

or all~\(m_{i} \geq \mathcal{l}\), in which case

\[r_{\mathcal{l}}\left( n \right) < 1.\]

Equation (85) then proves (83).~

\textbf{Corollary 8. Exact second moment and variance}

Assume~\(n \geq 2\). Then

\[\begin{matrix}
 & \mathbb{E}\left\lbrack x_{T}^{2} \right\rbrack = \mu_{T} + \frac{\prod_{i = 1}^{T}m_{i}\left( m_{i} - 1 \right)}{\left\lbrack n\left( n - 1 \right) \right\rbrack^{T - 1}}. & & \mathrm{(86)} \\
\end{matrix}\]

Consequently,

\[\begin{matrix}
 & \text{Var}\left( x_{T} \right) = \mu_{T} + \frac{\prod_{i = 1}^{T}m_{i}\left( m_{i} - 1 \right)}{\left\lbrack n\left( n - 1 \right) \right\rbrack^{T - 1}} - \mu_{T}^{2}. & & \mathrm{(87)} \\
\end{matrix}\]

If~\(\mu_{T} > 0\), this can equivalently be written as

\[\begin{matrix}
 & \text{Var}\left( x_{T} \right) = \mu_{T}\left\lbrack 1 + \frac{\prod_{i = 1}^{T}(\left. \ m_{i} - 1 \right)}{\left( n - 1 \right)^{T - 1}} - \frac{\prod_{i = 1}^{T}m_{i}}{n^{T - 1}} \right\rbrack. & & \mathrm{(88)} \\
\end{matrix}\]

\textbf{Proof.}~The identity

\[x_{T}^{2} = x_{T} + \left( x_{T} \right)_{2}\]

gives

\[\mathbb{E}\left\lbrack x_{T}^{2} \right\rbrack = \mu_{T} + F_{2}\left( T \right).\]

Taking~\(\mathcal{l} = 2\)~in (81) yields

\[F_{2}\left( T \right) = \frac{\prod_{i = 1}^{T}m_{i}\left( m_{i} - 1 \right)}{\left\lbrack n\left( n - 1 \right) \right\rbrack^{T - 1}}.\]

This proves (86). Subtracting~\(\mu_{T}^{2}\)~gives (87), and factoring
out~\(\mu_{T}\)~gives (88).~

\textbf{Corollary 9. Variance under-dispersion and the sparse-overlap
regime}

Assume~\(n \geq 2\). If~\(\mu_{T} > 0\), then

\[\begin{matrix}
 & \text{Var}\left( x_{T} \right)\mathbb{< E}\left\lbrack x_{T} \right\rbrack. & & \mathrm{(89)} \\
\end{matrix}\]

More precisely,

\[\begin{matrix}
 & 0 \leq \mathbb{E}\left\lbrack x_{T} \right\rbrack - \text{Var}\left( x_{T} \right)\mathbb{\leq E}\left\lbrack x_{T} \right\rbrack^{2}. & & \mathrm{(90)} \\
\end{matrix}\]

Consequently, along any sequence of MAO models satisfying

\[\mu_{T}\mathbb{= E}\left\lbrack x_{T} \right\rbrack \longrightarrow 0,\]

we have

\[\begin{matrix}
 & \mathbb{E}\left\lbrack x_{T} \right\rbrack - \text{Var}\left( x_{T} \right) = o\left( \mathbb{E}\left\lbrack x_{T} \right\rbrack \right). & & \mathrm{(91)} \\
\end{matrix}\]

If~\(\mu_{T} = 0\), then~\(x_{T} = 0\)~almost surely.

\textbf{Proof.}~From

\[\text{Var}\left( x_{T} \right) = \mu_{T} + F_{2}\left( T \right) - \mu_{T}^{2}\]

and Theorem 6, we have

\[0 \leq F_{2}\left( T \right) < \mu_{T}^{2}\]

whenever~\(\mu_{T} > 0\). This proves (89). Furthermore,

\[\mu_{T} - \text{Var}\left( x_{T} \right) = \mu_{T}^{2} - F_{2}\left( T \right),\]

so

\[0 \leq \mu_{T} - \text{Var}\left( x_{T} \right) \leq \mu_{T}^{2},\]

which proves (90). Along the nontrivial terms of a sequence for
which~\(\mu_{T} > 0\), division by~\(\mu_{T}\)~gives

\[0 \leq \frac{\mu_{T} - \text{Var}\left( x_{T} \right)}{\mu_{T}} \leq \mu_{T} \longrightarrow 0.\]

This proves (91). If~\(\mu_{T} = 0\), then at least one sampled subset
has size zero, so their full intersection is empty almost
surely.~\(\square\)

\textbf{Remark.}~Theorem 6 establishes the following MAO
factorial-moment inequality for the fully occupied count:

\[\begin{matrix}
 & F_{\mathcal{l}}\left( T \right)\mathbb{= E}\left\lbrack \left( x_{T} \right)_{\mathcal{l}} \right\rbrack \leq \left\{ \mathbb{E}\left\lbrack x_{T} \right\rbrack \right\}^{\mathcal{l}} = F_{1}\left( T \right)^{\mathcal{l}},\ 1 \leq \mathcal{l} \leq n. & & \mathrm{(}\mathrm{92}\mathrm{)} \\
\end{matrix}\]

For general occupancy sets~\(B_{1},\ldots,B_{\mathcal{l}}\), an
inequality of the form

\[\begin{matrix}
 & F_{\mathcal{l}}\left( B_{1},\ldots,B_{\mathcal{l}} \right) \leq \prod_{j = 1}^{\mathcal{l}}F_{1}\left( B_{j} \right) & & \mathrm{(}\mathrm{93}\mathrm{)} \\
\end{matrix}\]

does not hold without additional restrictions. A natural question is
therefore to determine the conditions
on~\(B_{1},\ldots,B_{\mathcal{l}}\)~under which such an inequality is
valid. One possible starting point is the common-set
case~\(B_{1} = \cdots = B_{\mathcal{l}} = B\), particularly
when~\(B = \left\{ t \right\}\)~or~\(B = \{ t,\ldots,T\}\). Whether the
resulting inequality

\[F_{\mathcal{l}}\left( B,\ldots,B \right) \leq F_{1}\left( B \right)^{\mathcal{l}}\]

holds, and under what additional conditions, is left for separate
investigation.

\textbf{8. Computational Implementation and Validation}

The accompanying Python implementation uses exact rational arithmetic
and takes the generalized MAO representation in Theorem 5 as its basic
computational formula:

\[F_{\mathcal{l}}\left( B_{1},\ldots,B_{\mathcal{l}} \right) = \frac{G_{T}\left( B_{1},\ldots,B_{\mathcal{l}} \right)}{\left( n \right)_{\mathcal{l}}^{T - 1}}.\]

Exact-level and threshold factorial moments are obtained by
taking~\(B_{j} = \left\{ t \right\}\)~and~\(B_{j} = \{ t,\ldots,T\}\),
respectively. Raw moments are then recovered by Stirling transformation,
and means, cross moments, variances, and covariances are derived within
the same framework.

The transversal
sum~\(G_{T}\left( B_{1},\ldots,B_{\mathcal{l}} \right)\)~is evaluated by
a sparse dynamic program that accumulates the unnormalized MAO weights
over the~\(T\)~subsets and the at most~\(2^{\mathcal{l}}\)~possible
membership states of~\(\mathcal{l}\)~specified distinct population
elements. This recursion is algebraically equivalent to sequential
coefficient accumulation in the joint PGF, but it directly evaluates the
MAO representation without constructing symbolic polynomial objects. For
fixed~\(\mathcal{l}\), a direct implementation requires at
most~\(O\left( \left( T + 1 \right)^{\mathcal{l}} \right)\)~occupancy-state
storage
and~\(O\left( T2^{\mathcal{l}}\left( T + 1 \right)^{\mathcal{l}} \right)\)~arithmetic
updates. It is therefore practical for moderate fixed orders,
including~\(\mathcal{l} \leq 4\), as required for the first four raw
moments.

The implementation was validated in three complementary ways. First, for
tractable parameter settings, such as
\(n = 9\)~and~\(\left( m_{1},m_{2},m_{3},m_{4} \right) = \left( 3,2,2,1 \right)\),
all~\(979,776\)~ordered subset configurations were exhaustively
enumerated. For every~\(t = 1,\ldots,4\), the MAO and enumeration
results agreed exactly, as rational numbers, for the means, variances,
and second through fourth raw moments of
both~\(x_{t}\)~and~\(x_{\geq t}\). Second, for the equal-size
case~\(n = 12\),~\(T = 4\), and~\(m = 4\), the independently implemented
specialized formulas (21)--(28) and (41)--(47) agreed exactly with the
generalized MAO calculations for all occupancy
levels~\(0 \leq t \leq 4\), including joint probabilities, cross
moments, covariances, and variances.

Third, for the larger instance, such
as~\(n = 1000\)~and~\(\left( m_{1},m_{2},m_{3},m_{4} \right) = \left( 90,120,75,110 \right)\),
the exact MAO moments were compared with Monte Carlo estimates
from~\(30,000\)~independently generated configurations and were
consistent within simulation variability. The complete numerical output
is generated by the accompanying Python code, which is available in
Supplementary files.

\textbf{9. Conclusion}

This paper develops a unified exact moment theory for multi-set
allocation occupancy under independent fixed-size sampling from a finite
population. The framework covers the complete occupancy profile

\[\left( x_{0},x_{1},\ldots,x_{T} \right),\]

the exact-occupancy counts~\(x_{t}\), and the threshold
counts~\(x_{\geq t}\). Its fundamental object is the generalized MAO
factorial moment

\[F_{\mathcal{l}}\left( B_{1},\ldots,B_{\mathcal{l}} \right),\]

which gives a direct probabilistic interpretation to the underlying MAO
allocation quantities and provides a common basis for factorial moments,
raw moments, cross moments, variances, and covariances.

The one- and two-element PGFs

\[Q = \Phi^{\left( 1 \right)}\ \mathrm{\text{and}}\ H = \Phi^{\left( 2 \right)}\]

determine the mean vector and covariance structure of the occupancy
profile. The general joint
PGF~\(\Phi^{\left( \mathcal{l} \right)}\)~extends the construction to
arbitrary factorial-moment order~\(\mathcal{l}\). Consequently, the
moment formulas for~\(x_{t}\)~and~\(x_{\geq t}\)~are specializations of
the same finite-population construction, and the equal-size formulas
follow as explicit reductions of the general theory.

For the fully occupied count

\[x_{T} = \left. \mid\bigcap_{i = 1}^{T}A_{i} \right.\mid,\]

the factorial moments satisfy

\[E\left\lbrack \left( x_{T} \right)_{\mathcal{l}} \right\rbrack = F_{\mathcal{l}}\left( \left\{ T \right\},\ldots,\left\{ T \right\} \right) = \frac{\prod_{i = 1}^{T}(\left. \ m_{i}\left. \ ) \right.\ _{\mathcal{l}} \right.\ }{\left( n \right)_{\mathcal{l}}^{T - 1}} \leq \left( E\left\lbrack x_{T} \right\rbrack \right)^{\mathcal{l}}.\]

In particular,

\[\text{Var}\left( x_{T} \right) < E\left\lbrack x_{T} \right\rbrack\]

in the nondegenerate cases covered by the strict inequality. Moreover,
in the sparse-overlap regime in
which~\(E\left\lbrack x_{T} \right\rbrack \rightarrow 0\),

\[E\left\lbrack x_{T} \right\rbrack - \text{Var}\left( x_{T} \right) = o\left( E\left\lbrack x_{T} \right\rbrack \right),\]

and hence

\[\text{Var}\left( x_{T} \right) \sim E\left\lbrack x_{T} \right\rbrack.\]

Thus, although fixed-size sampling produces factorial under-dispersion
relative to the Poisson benchmark, the mean and variance are
asymptotically equivalent when the expected fully occupied intersection
becomes sparse.

The accompanying implementation verifies the computational formulas
through exhaustive enumeration, independent equal-size reductions, and
Monte Carlo simulation. These checks establish consistency among direct
finite counting, the unified MAO factorial-moment construction, its
specialized closed forms, and stochastic simulation. The resulting
theory supplies exact finite-population benchmarks for recurrent
occupancy across multiple fixed-size subsets. Extensions to general
product inequalities
for~\(F_{\mathcal{l}}\left( B_{1},\ldots,B_{\mathcal{l}} \right)\),
their equality conditions, and their connections with finite-population
negative dependence remain subjects for further study.

\textbf{Funding:}~This research was supported by the National Natural
Science Foundation of China (82473137, 82273978)

\textbf{Conflict of interest:}~The authors declare no competing
interests.

\textbf{Data availability:}~No empirical dataset was generated or
analyzed for the theoretical results presented in this article.

\textbf{Code availability:}~The python code used to validate the
formulas is provide in the supplementary material.

\textbf{Author contributions:}~XGM developed the hypothesis, performed
mathematical induction, wrote the manuscript. XYX discussed the details,
wrote and revised the manuscript.

\textbf{Acknowledgements}

An OpenAI language model (ChatGPT) was used for language editing and
limited assistance with the presentation of numerical results. The
author developed the mathematical theory and the core computational
implementation, independently verified the results, and takes full
responsibility for the content of the manuscript.

\textbf{Supplementary material}

The supplementary material contains Python code for reproducing the
computational examples and independently checking the moment formulas.
All mathematical definitions, statements, and proofs required for the
results are contained in the main text.

\textbf{References}

Alon, N., and Spencer, J. H. (2016). The Probabilistic Method, 4th
Edition. Wiley Series in Discrete Mathematics and Optimization.

Amand, J., Fehlmann, T., Backes, C., and Keller, A. (2019). DynaVenn:
web-based computation of the most significant overlap between ordered
sets. BMC Bioinformatics \emph{20}, 743.

Borsboom, D., Deserno, M. K., Rhemtulla, M., Epskamp, S., Fried, E. I.,
McNally, R. J., Robinaugh, D. J., Perugini, M., Dalege, J., Costantini,
G.\emph{, et al.} (2021). Network analysis of multivariate data in
psychological science. Nature Reviews Methods Primers \emph{1}, 58.

Cochran, W. G. (1991). Sampling Techniques, 3rd Edition. Wiley Series in
Probability and Statistics.

Dokuchaev, N. G. (2000). Equations for probability distributions of
local occupation time on a surface for diffusion processes and control
problems. Journal of Mathematical Sciences \emph{99}, 1075-1088.

Fortunato, S. (2010). Community detection in graphs. Phys Rep
\emph{486}, 75-174.

Hájek, J. (1960). Limiting Distributions in Simple Random Sampling from
a Finite Population. A MAGYAR TUDOMÁNYOS AKADÉMIA MATEMATIKAI KUTATÓ
INTÉZETÉNEK KÖZLEMÉNYEI \emph{5}, 361-374.

Hur, B., Kang, D., Lee, S., Moon, J. H., Lee, G., and Kim, S. (2019).
Venn-diaNet : venn diagram based network propagation analysis framework
for comparing multiple biological experiments. BMC Bioinformatics
\emph{20}, 667.

Lex, A., Gehlenborg, N., Strobelt, H., Vuillemot, R., and Pfister, H.
(2014). UpSet: Visualization of Intersecting Sets. IEEE Transactions on
Visualization and Computer Graphics \emph{20}, 1983-1992.

Mao, X. G. (2026). Limit Theory of the Multi-set Allocation Occupancy
(MAO) Distribution: Normal and Poisson Approximations via MAO Norm.
arXiv.

Mao, X. G., and Xue, X. Y. (2018). General hypergeometric distribution:
A basic statistical distribution for the number of overlapped elements
in multiple subsets drawn from a finite population. arXiv:180806924.

Mao, X. G., and Xue, X. Y. (2022). Mean, Variance and Asymptotic
Property for General Hypergeometric Distribution. arXiv.

Mao, X. G., and Xue, X. Y. (2025). The Multi-set Allocation Occupancy
function and inequality (MAO function and MAO inequality): the
foundation of Generalized hypergeometric distribution theory. arXiv.

Pavlopoulos, G. A., Kontou, P. I., Pavlopoulou, A., Bouyioukos, C.,
Markou, E., and Bagos, P. G. (2018). Bipartite graphs in systems biology
and medicine: a survey of methods and applications. GigaScience
\emph{7}, giy014.

Särndal, C.-E., Swensson, B., and Wretman, J. (1992). Model Assisted
Survey Sampling. Springer Series in Statistics, Springer New York, NY.

Vatutin, V. A., and Mikhajlov, V. G. (1982). Limit theorems for the
number of empty cells in an equiprobable scheme for group allocation of
particles. Theory of Probability \& Its Applications \emph{27}, 734-743.

Wilson, K. R., and Anderson, D. R. (1985). Evaluation of a Density
Estimator Based on a Trapping Web and Distance Sampling Theory. Ecology
\emph{66}, 1185-1194.

\end{document}